\documentclass[a4paper,oneside]{amsart}
\usepackage{amssymb, color, graphics, colordvi,euscript}
\usepackage{amscd}
\usepackage{graphicx}
\usepackage{epsfig}

\newcommand{\bb}{\mathbb}
\newcommand{\cc}{\bb C}
\newcommand{\rr}{\bb R}

\newcommand{\ff}{\bb F}
\newcommand{\z}{\bb Z}
\newcommand{\zz}{\z/2}
\newcommand{\pp}{\bb P}

\newcommand{\xsig}{(X,\sigma)}
\newcommand{\xr}{{X}(\rr)}

\newcommand{\xc}{{X}}

\newcommand{\pic}{\operatorname{Pic}}

\newcommand{\rk}{\operatorname{rk}}

\newcommand{\ox}{\mathcal O_X}
\newcommand{\kx}{\mathcal K_X}

\newcommand{\chiox}{\operatorname{\chi(\ox)}}

\newcommand{\refbooktitle}[1]{\textit{#1}}
\newcommand{\refjourtitle}[1]{\textit{#1}}
\newcommand{\refpapertitle}{}

\numberwithin{equation}{section}%
\newtheorem{theo}[equation]{Theorem}%
\newtheorem{prop}[equation]{Proposition}
\newtheorem{lem}[equation]{Lemma}

\theoremstyle{remark}
\newtheorem{rem}[equation]{Remark}

\newtheorem{exemple}[equation]{Example}

\newenvironment{demo}[1]{\noindent\textit{#1.\ }}{\qed\par}

\begin{document}

\title[Real elliptic surfaces]{Topological types of real regular\\ jacobian elliptic surfaces}

\author{Fr\'ed\'eric Bihan \and Fr\'ed\'eric Mangolte } 
\address{Laboratoire de Math\'ematiques,
Universit\'e de Savoie, 73376 Le Bourget du Lac Cedex, France, Fax: +33 (0)4 79 75 81 42}
\email{Frederic.Bihan@univ-savoie.fr, Frederic.Mangolte@univ-savoie.fr}

%%%%%%%%%%%%%%%%%%%%%%%%%%%%%%%%%%%%%%%%
\begin{abstract} We present the topological classification of real parts of real regular elliptic surfaces with a real section.
\end{abstract}
%%%%%%%%%%%%%%%%%%%%%%%%%%%%%%%%%%%%%%%%

\maketitle

%\begin{center}
%Email: mangolte@univ-savoie.fr\quad
%Fax: +33 (0)4 79 75 81 42
%\end{center}
%\bigskip

\begin{quote}\small
\textit{MSC 2000:} 14P25 14J27
\par\medskip\noindent
\textit{Keywords:} elliptic surface, real algebraic variety, trigonal curve
\end{quote}

%%%%%%%%%%%%%%%%%%%%%%%%%%%%%%%%%%%%%%%%
\section{Introduction}\label{sec:intro}
%%%%%%%%%%%%%%%%%%%%%%%%%%%%%%%%%%%%%%%%

By a real algebraic variety $X$, we mean a pair $(X,\sigma)$ where $X$ is a complex algebraic variety and $\sigma$ an
antiholomorphic involution acting on $X$. The involution $\sigma$ is often called a real structure and the fixed point set $X^\sigma$ is called the real part of $X$ and is denoted by $\xr$.
From now on, a surface will be assumed to be a nonsingular surface and a
{\it topological type} of real surfaces will be a class of surfaces with homeomorphic real parts.

Our main purpose is to make a step towards the Enriques-Kodaira classification of real algebraic surfaces.
More precisely, we are interested in the classification of the topological types of real algebraic surfaces
in a given complex deformation family. The case of real rational surfaces goes back to Comessatti, \cite{Co14}
as well as the case of real abelian varieties \cite{Co32}, (see also \cite[Chap. III, IV]{Si89}). 
For real ruled surfaces, see \cite[Chap. V]{Si89}. The case of K3 surfaces is due to Nikulin and Kharlamov \cite{Kh76,Ni80}.
More recently, the classification of topological types of real Enriques surfaces has been obtained by Degtyarev and Kharlamov
\cite{DK96}, while the classification for real hyperelliptic surfaces has been achieved by Catanese and Frediani \cite{CF03}.

The remaining class of real algebraic surfaces of special type is made
by the so-called properly elliptic surfaces. In this direction there is the fundamental
classification of the real singular fibres of an elliptic fibration by R.~Silhol \cite{Si84}
and some partial results about the global classification in \cite{Ma00} and \cite{AM06}.

A {\em real elliptic surface} will be a morphism $\pi \colon X\to \pp^1$ defined over $\rr$,
where $X$ is a real algebraic surface such that for all but finitely many points $u \in  \pp^1$,
the fibre $X_u= \pi^{-1}(u)$ is a nonsingular curve of genus one. When $\pi$ admits at least one singular fibre,
the surface $X$ is {\em regular}, which means that $H^1(X,\ox)=\{0\}$. 
The elliptic fibration $\pi \colon X\longrightarrow \pp^1$ will be called {\em relatively  minimal} if no fibre of
$\pi$ contains an exceptional curve of the first kind. 
When a (real) relatively minimal elliptic surface $\pi \colon X\to \pp^1$ admits a (real) section $s \colon \pp^1 \to X$ we call $X$ a
{\em (real) jacobian elliptic surface}. 
\medskip

Over $\cc$, two regular relatively minimal elliptic surfaces with no multiple fibres are deformation equivalent if and only if
their holomorphic Euler characteristics are equal, see~\cite{Kas77}.
  
The main goal of this paper is to present the classification of the topological types of real regular jacobian elliptic surfaces
in each complex deformation family of regular jacobian elliptic surfaces. 

In order to describe the topological types, we denote by $S_g$
the smooth orientable surface of genus $g$, by $S$ the two dimensional sphere ($S=S_0$) and by
$V_q$ the non orientable surface of Euler characteristic $2-q$ ($V_q$ is diffeomorphic
to the connected sum of $q$ copies of $\pp^2(\rr)$).
In what follows, the disjoint union of two surfaces $A$ and $B$ is denoted by $A \sqcup B$ while the disjoint union
of $a$ homeomorphic copies of $A$ is denoted by $a \, A$.

An abstract {\em Morse simplification} on a topological type
is a transformation which decreases the total Betti number by 2. There are
two kinds of Morse simplifications:
\begin{itemize}
\item removing one spherical component, $S \to \emptyset$,
\item contracting one handle $S_{g+1} \to S_g$ or $V_{q+2} \to V_q$.
\end{itemize}

A topological type of a real regular jacobian elliptic surface is called {\em extremal} if it cannot be obtained
by a Morse simplification from the topological type of another real regular jacobian elliptic surface in the
same complex deformation family (i.e. with the same holomorphic Euler characteristic).

\begin{theo}\label{theo:classif}
Let $k \geq 1$ be an integer. The extremal topological types of real regular jacobian elliptic surfaces
with holomorphic Euler characteristic $k$ are:

\begin{enumerate}
	\item $M$-surfaces, 
	$a = k + 4\lambda - 1$, $l = 5k - 4\lambda$, $\lambda = 0,1,\dots , k$
	\begin{itemize}
		\item $S_l \sqcup a\, S$,  $k$ even or
		\item $V_{2l} \sqcup aS$, $k$ odd.
	\end{itemize}\medskip
	\item $(M-2)$-surfaces,
	$a = k + 4\lambda$, $l = 5k - 4\lambda -3$, $\lambda = 0,1,\dots , k - 1$
	\begin{itemize}
		\item $S_l \sqcup a\, S$,  $k$ even or
		\item $V_{2l} \sqcup aS$, $k$ odd.
	\end{itemize}\medskip

	\item
	\begin{itemize}
		\item $S_1 \sqcup S_1$,  $k$ even or
		\item $V_2 \sqcup V_2$, $k$ odd.
	\end{itemize}
\end{enumerate}

Any topological type of a real regular jacobian elliptic surface with holomorphic Euler characteristic $k$
is the result of a sequence of Morse simplifications from one of
the previous extremal types.

Conversely, any topological type corresponding to a surface with first Betti number $\geq 2$ and
which can be obtained by a sequence of Morse simplifications from one of the previous extremal topological
types is the topological type of a real regular jacobian elliptic surface $X$ with $k = \chiox$.
\end{theo}

Note that in the definition of a jacobian surface, we supposed
that the elliptic fibration was relatively minimal. In fact, the analogue of Theorem~\ref{theo:classif}
without this hypothesis, that is, the classification of topological types
of real regular elliptic surfaces with a real section and contained in a given complex deformation family,
can be deduced directly from Theorem~\ref{theo:classif}.
Indeed, over $\cc$, two regular elliptic surfaces with no  
multiple fibres are deformation equivalent if and only if their
holomorphic Euler characteristics are equal {\em and} their canonical classes have  
the same degree. To realize a topological type in  one complex family, say $k = \chiox$ and $K_X^2 = -m < 0$, take any real regular jacobian elliptic surface $Y$ with holomorphic  Euler characteristic $k$, then by definition, $K_Y^2 = 0$. Let $X$ be the blow-up of $m$ points which are globally fixed  
by the real structure. Then $K_X = -m$. Each blow-up  
centered in a real point gives rise to a connected sum with a real  
projective plane. Conversely, any topological type of a real regular  
elliptic surface with a real section can be obtained in this way.

The rest of the paper is devoted to the proof of Theorem~\ref{theo:classif} (see Subsection~\ref{subsection:finalproof}).
\medskip

Let $\pi \colon X\longrightarrow \pp^1$ be a real regular jacobian elliptic surface. Denote by $k = \chiox$ the holomorphic Euler characteristic of $X$.
The Betti numbers  of $\xr$ are subject to prohibitions. Classical inequalities and congruences on
real regular algebraic surfaces (see~\cite{DK00} for a survey) yield a finite list for the allowed pairs of topological Euler characteristic $\chi(\xr)$
and total Betti number $h_*(\xr)$ with coefficients in $\zz$, see Section~\ref{sec:restrict}. The resulting list is shortened thanks to the following inequalities
which arise essentially because the surface
is elliptic. Namely the first  mod 2 Betti number $h_1(\xr) = \rk \left( H_1(\xr , \zz) \right)$ and the number of connected components are bounded by 
\begin{equation}\label{ineq:regular}
h_1(\xr) \leq 10 k\;, \qquad
\#\xr \leq 5 k\;.
\end{equation}

The upper bound on the first Betti number is due to a more general result of V.~Kharlamov (see Theorem~\ref{theo:rv}).
The upper bound on the number of connected components is due to M.~Akriche and the second author, see \cite{AM06}. 

Using the general prohibitions on real algebraic surfaces and the two specific inequalities (\ref{ineq:regular}), we draw for each $k$,
the diagram of the allowed pairs $(\chi  , h_*)$, see Figure~ \ref{fig:diammond}. 
The existence of a real section imposes strong conditions on the real part of $X$.
Namely, $\xr \ne \emptyset$, one connected component of $\xr$ has Euler characteristic $\leq 0$ and other
components are diffeomorphic to spheres, unless $\xr$ is the disjoint union of two tori (if $k$ even) or two Klein bottles (if $k$ odd), see Proposition~\ref{prop:jacotype}.
Hence, except for the pair $(0,8)$, there is a one-to-one correspondence between allowed pairs $(\chi  , h_*)$ and allowed topological types.
The list given in Theorem~\ref{theo:classif} is exactly the list of topological types obtained from the diagram in Figure~ \ref{fig:diammond}
using this correspondence.

Conversely, for each $k \geq 1$ and for each topological type listed in Theorem~\ref{theo:classif}, we prove the existence of a real regular jacobian elliptic surface
with the desired topological type. These surfaces are obtained in the following way. 

Fix an integer $k \geq 1$ and let $R = \ff_{2k}$ be the Hirzebruch surface of degree $2k$ endowed with the real ruling
$\tilde \pi \colon R \to \pp^1$. Let $E$ be the unique nonsingular algebraic curve of $R$ with nonzero negative selfintersection.
By unicity, $E$ is defined over $\rr$. Consider a nonsingular real algebraic curve $C$ of bidegree $(3,0)$ on $R$
(see Section~\ref{sec:hirze} for the definition of bidegree). Such curves are often called {\it trigonal curves}. Let
$\rho \colon X \rightarrow R$ be the double covering ramified over the union $C \cup E$ which is of {\em even} bidegree $(4 , -2k)$.
Standard calculation on double coverings yields  $\chiox = k$, and it is an easy exercise to check that $\pi = \rho \circ \tilde \pi$ is an elliptic fibration.
In fact, $X$ is a jacobian elliptic surface.

Since the curve $C$ is real, the surface $X$ becomes real with respect
to the two complex conjugations which are lifts of the one on $R$ (they are interchanged by the covering involution).
The topological type of the elliptic surface with respect to either complex conjugations is determined by the real scheme
of $C$, i.e. the pair $(F(\rr), C(\rr))$ up to homeomorphism, and the parity of $k$ (see Lemma~\ref{L:curvetosurf}). A Morse simplification for
either real elliptic surfaces corresponds then to collapsing one oval of the curve $C$.
Hence, the list in Theorem~\ref{theo:classif} results in a list of real schemes of real trigonal
curves on $R$, where extremal ones are those which cannot be obtained from another by collapsing one oval.

Trigonal curves which produce the exceptional topological types of two tori (if $k$ is even) or two Klein bottles
(if $k$ is odd) are easy to construct, see Example~\ref{E:specialextremal}. 
The other extremal cases are constructed
by means of the {\it combinatorial patchworking}. The combinatorial patchworking is a combinatorial version
of the {\it Viro method}, which is a powerful construction method of real algebraic varieties with prescribed topology, see ~\cite{Vi83, Vi84, Vi}.
For the remaining real schemes, we use the fact that ovals of trigonal curves can be collapsed independently
\footnote{It is worth noting that this is not true for other types of curves, e.g plane projective curves of a given degree}.
This can be shown using the theory of {\it dessins d'enfants}. Though this statement can be easily deduced from~\cite{Or},
we give here a detailed proof for the reader's convenience (see Proposition~\ref{P:collapsing}).
The emergence of dessins d'enfants in real algebraic geometry is recent~\cite{NSV,Or}
(see also, for example,~\cite{BeBr,Bi}). It is S. Orevkov~\cite{Or} who discovered that dessins d'enfants can be used to construct real trigonal curves
with given real scheme. Recently, B. Bertrand and E. Brugalle~\cite{BeBr} showed how one can recover the dessin d'enfant of a trigonal curve
constructed by the combinatorial patchworking.

One advantage of the construction used in Section~\ref{sec:collapsing}
comes from the obvious observation that a Morse simplification may not be obtained a priori in a continuous family of real regular
jacobian elliptic surfaces. This construction shows that such a family does exist in each
case except possibly for the exceptional extremal pair of tori or Klein bottles.

Let us finally point out the recent appearance of~\cite{DIK06} in which dessins d'enfants are used in order to study real elliptic surfaces
up to equivariant deformation.

\bigskip
We want to thank E. Brugall\'e, F. Catanese, A. Degtyarev, I. Itenberg and V. Kharlamov for fruitful discussions.
The second author thanks also F. Catanese for a stimulating visit in Bayreuth.

%%%%%%%%%%%%%%%%%%%%%%%%%%%%%%%%%%%
\section{Restrictions on real regular elliptic surfaces}\label{sec:restrict}
%%%%%%%%%%%%%%%%%%%%%%%%%%%%%%%%%%%

In this section, we collect some classical restrictions on the topology of a real algebraic surface $\xsig$ in term of the numerical invariants of the complex surface $X$. Next, we apply them to a real regular elliptic surface. For a complex regular elliptic surface $X$, most of the numerical invariants can be deduced from the holomorphic Euler characteristic  $\chiox$. Recall that in fact, when $X$ is jacobian, its complex deformation class is determined by $\chiox$, see~\cite{Kas77}.

\medskip
Let $\xsig$ be a real algebraic manifold. Denote by $B_i(X) = \rk \left( H^i(\xc , \zz \right)$ the $i$-th mod 2 Betti number of $\xc$ and by 
$$
	h_i(\xr) = \rk \left( H_i(\xr , \zz) \right)
$$
the $i$-th mod 2 Betti number of $\xr$. Let $h_*(\xr) = \sum h_i(\xr)$ and $B_*(X) = \sum B_i(X)$.
The first important prohibition is given by the Milnor-Smith-Thom inequality.

\begin{prop}\label{prop:smith}\label{eq:smith}
Let $\xsig$ be a real algebraic manifold. Then
	$$
		h_*(\xr) \leq B_*(X) \;;
	$$
	$$
	B_*(X) - h_*(\xr) \, \textrm{is even.}
	$$
\end{prop}

\begin{demo}{Proof}
See e.g. \cite[Chap. II]{Si89} or \cite[3.1.1]{DK00}. 
\end{demo}

\medskip
Let us denote by $\tau(X)$ the signature of the cup-product form on the real vector space $H^2(\xc , \rr)$ and by $\tau^-(X)  = \frac 12(B_2(X) - \tau(X))$ its negative index.
The following inequality, often called Comessatti's inequality, is well-known.

\begin{prop}\label{prop:comessatti}
Let $\xsig$ be a real algebraic surface, then 
$$
\vert \chi(\xr) - 1 \vert \leq \tau^-(X)
$$ 
\end{prop}

\begin{demo}{Proof}
See e.g. \cite[3.1.2]{DK00}.
\end{demo}

\medskip
An important class of restrictions on extremal real surfaces is given by the Rokhlin congruences and their generalizations by Gudkov and Kharlamov. Recall that $B_*(X) - h_*(\xr)$ is even and let 
$2d = B_*(X) - h_*(\xr)$. A real algebraic surface is an {\em $M$-surface} if $d = 0$ and an {\em $(M-d)$-surface} if $d \ne 0$.

\begin{prop}\label{prop:cong}
	If $\xsig$  is an $M$-surface $($that is $d = 0)$, then
	$$
		\chi(\xr) \equiv \tau(X) \mod 16\;.
	$$

	When $\xsig$  is an $(M -1)$-surface $($that is $d = 1)$, then
	$$
		\chi(\xr) \equiv \tau(X) \pm 2 \mod 16\;.
	$$
\end{prop}

\begin{demo}{Proof}
See e.g. \cite[Chap. II]{Si89} or \cite[2.7.1]{DK00}.
\end{demo}

\begin{prop}\label{prop:congellipt}
	Let $\pi \colon X\to \pp^1$ be a real regular relatively minimal elliptic surface with no multiple fibres
	
When $\xsig$  is an $M$-surface, then
$$
	\chi(\xr) \equiv 8\chiox \mod 16\;.
$$

When $\xsig$  is an $(M -1)$-surface, then
$$
	\chi(\xr) \equiv 8\chiox \pm 2 \mod 16\;.
$$
\end{prop}

\begin{demo}{Proof} The elliptic fibration has no multiple fibres thus the homology of the complex surface $X$ has no torsion. Whence we will not distinguish between the Betti numbers of $X$ and the mod 2 Betti numbers of $X$. As $\pi$ is relatively minimal, the first Chern number $c_1^2(X)$ is zero and the Noether's formula yields
\begin{equation}\label{eq:noether}
	\chi_{top}(\xc) = 12\chiox\;.
\end{equation}
For $X$ a regular algebraic surface (i.e. with $q(X) = \dim H^1(X,\ox) = 0$), the signature of the cup-product form is 
$\tau(X) = 4\chiox - 2 -B_2(X)$,
see e.g. \cite[Thm. 2.7]{BPV}.
Moreover, a regular surface has Betti numbers $B_1(X) = B_3(X)= 0$  whence
$\tau(X) = -8\chiox$. The conclusion follows from Proposition~\ref{prop:cong}.
\end{demo}

\begin{rem}
Let $k=\chiox$. From the previous congruences, the Comessatti inequality \ref{prop:comessatti}, and the Theorems~\ref{theo:rv} and \ref{theo:cc} below, we get the list of the $k + 1$ admissible pairs of Betti numbers for $M$-surfaces: 
$$
(\#\xr , h_1(\xr)) = (k +4\lambda, 10k -8\lambda)\;, \quad\lambda \in \{0,1,\dots,k\}\;.
$$

From the Betti numbers, we will deduce the jacobian topological type by Proposition~\ref{prop:jacotype}.
\end{rem}

\bigskip
By Morse simplifications, any $(M-1)$-jacobian topological type allowed by the Proposition \ref{prop:congellipt} comes from an $M$-jacobian topological type, otherwise said, there are no extremal  $(M-1)$-jacobian topological type. There are $k$ extremal $(M-2)$-jacobian topological type, see Theorem~\ref{theo:classif}. To prove that the extremal topological types given in Theorem~ \ref{theo:classif} are the only ones, the former classical restrictions on real algebraic surfaces applied to a real regular relatively minimal elliptic surface with no multiple fibres are not sufficient. We need the following additional restrictions.

\begin{theo}\cite{AM06,Ma00}\label{theo:rv}\label{theo:cc}
	Let $X \to \pp^1$ be a real regular relatively minimal elliptic surface with no multiple fibres, 	then
	$$
		h_1(\xr) \leq 10\chiox \;.
	$$

	Moreover, if $\pi$ is jacobian, then
	$$
		\#\xr \leq 5\chiox \;.
	$$
\end{theo}

\begin{prop}\label{prop:ineqellipt}
	Let $\pi \colon X\to \pp^1$ be a real regular jacobian elliptic surface. We have 
	\begin{eqnarray*}
		&4 \leq h_*(\xr) \leq 12\chiox\\
		&2 - 10\chiox  \leq \chi(\xr) \leq 10\chiox - 2
	\end{eqnarray*}
	
\end{prop}

\begin{demo}{Proof}
The inequality (\ref{eq:smith}) and the Noether formula for a relatively minimal elliptic fibration (\ref{eq:noether}) yield
$$
	h_*(\xr) \leq 12\chiox \;.
$$

Furthermore, there is a connected component $N_0$ of $\xr$ containing the real part $S(\rr)$ of the image of the real section $s \colon \pp^1 \to X$. Let $\alpha \in H_1(N_0, \zz)$ be the fundamental class of $S(\rr)$ and $\beta \in H_1(N_0, \zz)$ the class of the restriction to $N_0$ of a smooth fibre, by hypothesis $\alpha \cap \beta = 1 \mod 2$  whence $h_1(N_0) \geq 2$.

From $B_*(X) = 12 \chiox$,  a regular relatively minimal elliptic surface with no multiple fibres $X$ has $B_2(X) = 12 \chiox - 2$ and $\tau^-(X) = 10 \chiox - 1$. We get from Proposition~\ref{prop:comessatti}, 
$$
	2 - 10\chiox  \leq \chi(\xr) \leq 10\chiox \;.
$$

We want to prove that, in fact, $\chi(\xr) \leq 10\chiox - 2$.
Suppose for a moment that $\chi(\xr) = 10 \chiox$, that is $\#\xr = 5\chiox + \frac 12 h_1(\xr)$. Now from the Theorem~\ref{theo:rv}, we get the inequality $\#\xr \leq 5 \chiox$ thus $\#\xr = 5 \chiox$ and $h_1(\xr) = 0$.  
But we have seen previously that $\xr$ contains a component $N_0$ with $h_1(N_0) \geq 2$; a contradiction. 

For $X$ a real relatively minimal elliptic surface, the Euler characteristic of $\xr$ is even. Indeed, by  Proposition~\ref{prop:smith}, $\chi(\xr) \equiv \chi_{top}(\xc) \mod 2$ and by (\ref{eq:noether})  $\chi_{top}(\xc)$ is even. Thus  $\chi(\xr) = 10 \chiox-1$ does not occur and finally this yields $\chi(\xr) \leq 10\chiox - 2$.
\end{demo}
\medskip

The previous proof shows in particular that if $X$ is a real regular jacobian elliptic surface,
then $\xr \neq \emptyset$ and $h_1(\xr) \geq 2$. This yields

\begin{equation*}\label{E: lowerbounds}
4-\chi(\xr) \leq h_*(\xr) \quad \mbox{and} \quad \chi(\xr) \leq h_*(\xr) - 4.
\end{equation*}

We also know that $\chi(\xr)$ and $h_*(\xr)$ are even integer numbers.
Together with Propositions~\ref{prop:congellipt} and \ref{prop:ineqellipt} and Theorem~\ref{theo:rv},
the preceding analysis allows us to draw, for each $k$, the diagram of possible values for the pair $(\chi(\xr),h_*(\xr)$ when
$X$ is a real regular jacobian elliptic surface.

\begin{figure}[htbp]
\bigskip

\begin{center}
\includegraphics{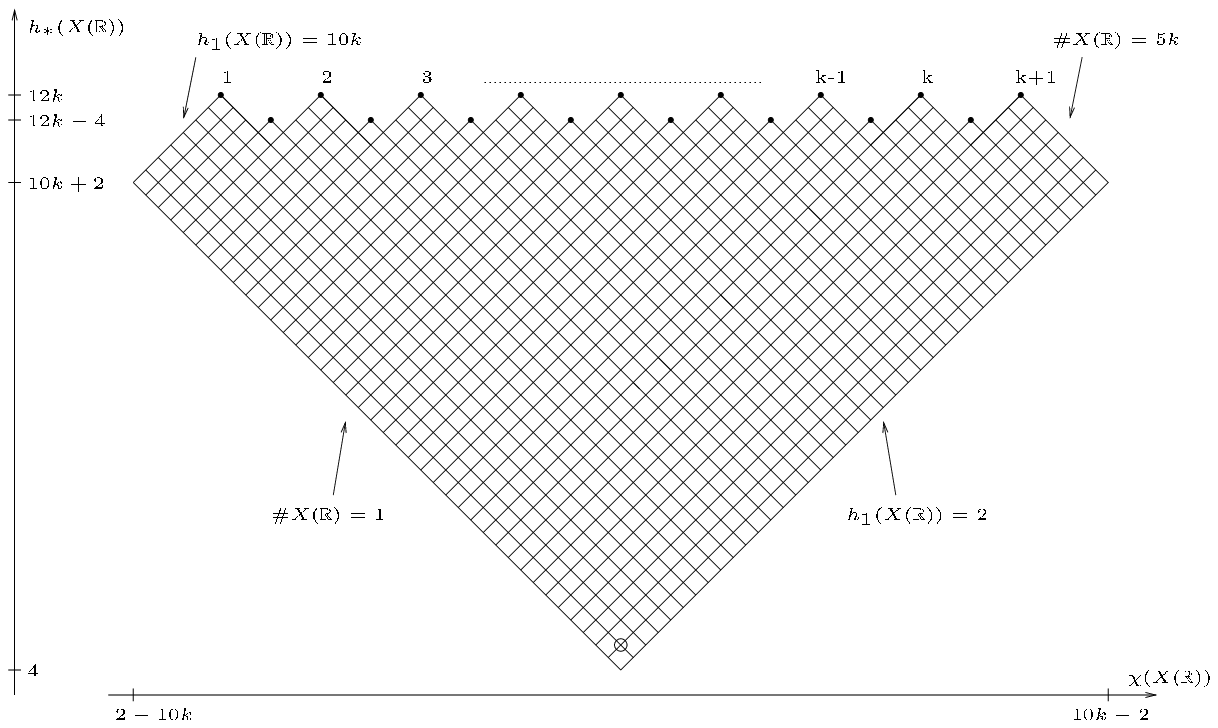}
 \caption{Diagram of possible values for $(\chi(\xr),h_*(\xr)$}\label{fig:diammond}
	\end{center}
	\raggedright{
{\scriptsize
The diagram is the  right one for $k=8$. The black points correspond to the $k+1$ $M$-surfaces (which are automatically extremal)
and the extremal $(M-2)$-surfaces. The small circle surrounds the point $(0,8)$ which is the only point
in the diagram corresponding to two possible topological types comprising an extremal one: $S_1 \sqcup S_1$ if $k$ is even
or $V_2 \sqcup V_2$ if $k$ is odd. Note that for any $k$ the diagram is symmetric
with respect to the vertical axis $\chi(\xr)=0$ while this axis contains a point corresponding to an $M$-surface (as it is the case in the picture) only when $k$ is even.}}
	\label{f:coloring}
\bigskip

\end{figure}

From the diagram, we will infer the list of allowable topological types of real regular jacobian elliptic surfaces. For that we use the following result.

\begin{prop}\label{prop:jacotype}
Let $\pi \colon X\to \pp^1$ be a real regular jacobian elliptic surface. Let $k=\chi(\ox)$.
Then $\xr$ is non empty and it is orientable if $k$ is even or non orientable if $k$ is odd. Moreover, $\xr$ has one connected component of Euler characteristic
$\leq 0$ and other components are diffeomorphic to a sphere, unless $\xr$ is the disjoint union of two tori (if $k$ even) or two Klein bottles (if $k$ odd).
\end{prop}

\begin{demo}{Proof}
From \cite[Chap. VII]{Si89}, we get that when the elliptic fibration admits a real section, any real singular fibre has connected real part except
if the fibre is real isomorphic to the real cubic curve with one isolated node. Recalling that the real part of a smooth real curve of genus one has
at most two connected components, the result easily follows except for the assertion about orientability which is deduced from the following Lemma.
\end{demo}

\begin{lem}\label{lem:ornor}
	Let $X \to \pp^1$ be a real regular relatively minimal elliptic surface with no multiple fibres.

	If $\chi(\ox)$ is even, $\xc$ is a $\operatorname{Spin}$-surface, and if $\xr\ne\emptyset$,
	$$
		\xr \textrm{ is orientable.}
	$$

	If $\chi(\ox)$ is odd, $\xc$ is a non $\operatorname{Spin}$-surface, and if $X \to \pp^1$ admits a real section or if $\xsig$ is an $M$-surface,
	$$
		\xr\ne\emptyset \textrm{ and $\xr$ is non orientable.}
	$$
\end{lem}

\begin{demo}{Proof}
	By \cite[p. 162]{BPV}, we have for a canonical divisor of $X$~:
	$$
		\kx = (\chi(\ox)-2) \, \tilde{F}
	$$
	where $\tilde{F}$ is any fibre of $X \to \pp^1$. From $c_1(T\xc) = -c_1 (\kx)$ we obtain for the second Stiefel-Whitney class of the tangent bundle of $X$
	$$
		w_2(\xc) = \chiox \mod 2\;.
	$$
	We use \cite[5.7]{Ma00} for the first two assertions.
	We use \cite[3.4.1(1)]{DK00} for the last assertion on $M$-surfaces.
\end{demo}

\medskip

Proposition~\ref{prop:jacotype} implies that, with the only exception of the point $(0,8)$,
each point in the diagram depicted in Figure~\ref{fig:diammond} corresponds to at most one topological type
of a real regular jacobian elliptic surface.

Namely, if $k$ is even the point $(0,8)$ corresponds to  $S_2 \sqcup S$ (non extremal) and to $S_1 \sqcup S_1$ (extremal).
If $k$ is odd, the point $(0,8)$ corresponds to $V_4 \sqcup S$ (non extremal) and to $V_2 \sqcup V_2$ (extremal).
Any point $(\chi, h_*) \neq (0,8)$ of the diagram corresponds to
the topological type $S_l \sqcup aS$ if $k$ is even, or to $V_{2l} \sqcup aS$ if $k$ is odd, where $a$ and $l$ are determined by
$\chi=2(a+1)-2l$ and $h_*=2(a+1)+2l$.

Collecting the results of this section, we thus obtain a list of non prohibited
topological types, which coincides with the one given in Theorem~\ref{theo:classif}.

%%%%%%%%%%%%%%%%%%%%%%%%%%%%%%%%%
\section{Hirzebruch surfaces and real elliptic surfaces}\label{sec:hirze}
%%%%%%%%%%%%%%%%%%%%%%%%%%%%%%%%%

Let $n \geq 1$, denote by $\tilde{\pi} \colon \ff_n \rightarrow \pp^1$ the $\pp^1$-bundle over $\pp^1$ associated with the Hirzebruch surface $\ff_{n}$.
Let $(u,x)$ be affine coordinates on $\ff_{n}$ so that $\tilde{\pi}$ can be written locally as $(u,x) \mapsto u$.
Denote by $B$ and $F$ the classes in the second homology group $H_2(\ff_n, \z)$ of the image of a section disjoint from the exceptional divisor $E$ and a fibre of $\tilde{\pi}$, respectively. 
As it is well-known, the homology group $H_2(\ff_n, \z)$ is generated by
$B$ and $F$ with intersection numbers $B \cdot B = n$, $F \cdot F = 0$ and $B \cdot F = 1$.

A curve $C$ on $\ff_n$ is of bidegree $(a,b)$ if $[C] = aB + bF$ in $H_2(\ff_n, \z)$, or equivalently, if $a = [C] \cdot  F$ and $b = [C] \cdot E$.
The intersection number of two curves of bidegree $(a,b)$ and $(c,d)$ is then equal to $acn + ad +bc$. 
The exceptional divisor $E$ has bidegree $(1,-n)$ and self-intersection
$-n$. 

The Hirzebruch surface $\ff_n$ inherits a real structure coming from
the usual complex conjugation on some affine chart of $\ff_{n}$. 
This is the standard real structure
on $\ff_n$ and we will consider only this one. The real part $\ff_n(\rr)$ of $\ff_n$ is homeomorphic to a two-dimensional torus
if $n$ is even, or to a Klein bottle if $n$ is odd.

\medskip
A {\it trigonal curve} is a curve of bidegree $(3,0)$ on $\ff_n$. 
Let $r \in \rr[u,x]$ be a polynomial whose Newton polygon (the convex hull of the exponent vectors of the monomials appearing with non zero coefficients)
is the triangle with vertices $(0,0)$, $(3n,0)$ and $(0,3)$.
Such a polynomial defines (affinely) a trigonal curve $C$. A classical (non linear) elimination process allows to eliminate the monomials
$u^lx^2$, $l \in \{ 0,\dotsc,n \}$, from $r$ in order to get a {\it reduced} affine equation for $C$

\begin{equation}\label{E:reduced}
x^3+p(u)x+q(u)=0,
\end{equation}
where $p$ and $q$ are polynomials of degrees $2n$ and $3n$, respectively.
The discriminant of this polynomial, viewed as a polynomial in $\rr[u][x]$, is
$$
\Delta=4p^3+27q^2 \in \rr[u] \;.
$$
The roots of $\Delta$ are those $u_0$ for which the vertical line $u = u_0$ is either tangent to, or contains a singular point of, the curve $C$. We have $\Delta(u_0)<0$ exactly when the line $u = u_0$ intersects $C(\rr)$ in three points.

Consider the Hirzebruch surface $\ff_{2k}$ of even degree $2k$ and let $C$ be a real nonsingular trigonal curve on $\ff_{2k}$.
A connected component of $C(\rr)$ is called a {\it pseudo-line} if it intersects each vertical line $u=u_0$,
otherwise it is called an {\it oval}.
An oval separates $\ff_{2k}(\rr)$ into two connected components, one of them being homeomorphic to a disk called {\it interior}
of the oval. The image of a pseudo-line by the projection $\ff_{2k}(\rr) \rightarrow \pp^1(\rr)$ induced by $(u,x) \mapsto u$ is the whole
$\pp^1(\rr)$, while the image of an oval is an interval whose endpoints are real roots of the discriminant $\Delta$.

With one exception, all but one connected components of $C(\rr)$ are ovals, while the remaining one is a pseudo-line.
The exception arises when $C(\rr)$ consists of three (non intersecting) pseudo-lines. This happens only if $\Delta < 0$ on $\pp^1(\rr)$.
We give here an explicit example which will be used for the construction of extremal elliptic surfaces
with topological type $S_1 \sqcup S_1$ (if $k$ even) or $V_2 \sqcup V_2$ (if $k$ is odd), see Lemma~\ref{L:curvetosurf} below.

\begin{exemple}\label{E:specialextremal}
Consider the non reduced polynomial

$$
r(u,x) =  (x-g_1(u)) \cdot (x-g_2(u)) \cdot (x-g_3(u)),
$$
where $g_1, g_2$ and $g_3$ are real polynomials of degree $2k$
such that $g_2-g_1$, $g_3-g_1$ and $g_3-g_2$ have each $2k$ (simple) non real roots
and the resulting $6k$ roots are pairwise distinct. In particular, $g_2-g_1$, $g_3-g_1$ and $g_3-g_2$
are polynomials of degree $2k$ with no real roots. This implies that $r$ defines a trigonal curve whose real part
consists of three non intersecting pseudo-lines. This trigonal curve has as singularities $6k$
non real ordinary double points coming from the above $6k$ non real roots. A small perturbation
gives a nonsingular trigonal curve whose real part consists of three pseudo-lines.
\end{exemple}

\begin{figure}[htbp]
\bigskip

\begin{center}
\includegraphics{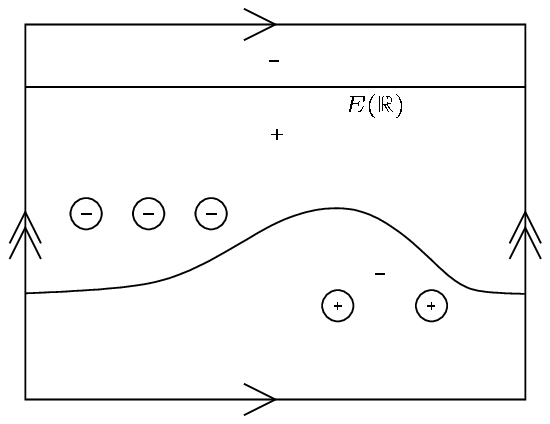}
\caption{Real scheme $\langle 2 \mid 3 \rangle$}
\label{F:realscheme}
\end{center}
\bigskip

\end{figure}

Assume now that $C(\rr)$ does not consist of three pseudo-lines. Let 
$$
H = 0
$$ 
be some bi-homogeneous polynomial equation for $C \cup E$.
We have $[C] \cdot E = 0$, hence $C$ does not intersect $E$ and $C \cup E$ is a real nonsingular curve of bidegree $(3,0)+(1,-2k)=(4,-2k)$.
Since this bidegree is even, the polynomial $H$ has a well-defined sign on $\ff_{2k}(\rr)$. Hence the ovals
of $C$ can be divided into two groups according to the sign of $H$ in the interior
of each oval. It follows that,
with the exception of $C(\rr)$ consisting of three pseudo-lines, the real scheme of $C$, that is the pair
$(\ff_{2k}(\rr), C(\rr))$ up to homeomorphism, is uniquely determined
by the numbers $a$ and $b$ of ovals in each of these two groups and will be encoded
by $\langle a \mid b \rangle$, see Figure~\ref{F:realscheme}. Note that these two numbers are defined only up to permutation.

Let $\rho \colon X \mapsto \ff_{2k}$ be the double covering of $\ff_{2k}$ ramified over $C \cup E$ (it is well-defined
since the bidegree of the ramification curve is even). The complex surface $X$ is a nonsingular jacobian elliptic surface, which can be equipped with two real structures lifting that of $\ff_{2k}$. According to the chosen real structure, the real part
of $X$ projects via $\rho$ onto the subset of $\ff_{2k}(\rr)$ where either $H \geq 0$, or $H \leq 0$.

\begin{lem}\label{L:curvetosurf}
If $C(\rr)$ consists of three pseudo-lines, then $\xr$ is homeomorphic to
$S_1 \sqcup S_1$ if $k$ even, or to $V_2 \sqcup V_2$ if $k$ is odd.

If $C$ has real scheme $\langle a \mid b \rangle$, then, up to permutation of $a$ and $b$ which corresponds to a permutation
between the two real structures for $X$, the surface $\xr$ is homeomeomorphic to
$a \, S \sqcup S_{b+1}$ if $k$ is even, or to $a \, S \sqcup V_{2b + 2}$ if $k$ is odd.
\end{lem}
\begin{proof}
A direct computation shows that $\chi(\xr)=0$ if $C(\rr)$ consists of three pseudo-lines
and $\chi(\xr) = 2(a-b)$ if $C$ has real scheme $\langle a \mid b \rangle$.
It suffices then to determine if each non spherical component is orientable or not. 
Let $M$ be such a component, and let $S$ be a real section whose real part $S(\rr)$ is contained in $M$. The linear class of the effective divisor $S$ is the
(reduced part) of the pull-back of the class $B$ in $\pic (\ff_{2k})$.  The class $B$ is a part of the branch class of the double covering $\rho$. Hence
$[S] \cdot [S] = \frac 12 B \cdot B = k$, see e.g. \cite{Pe81}. Thus, when $k$ is odd, the square $\alpha^2$ of the fundamental class of $S(\rr)$ in $H_1(M, \zz)$
does not vanish  whence $M$ is non orientable.

Conversely, suppose $k$ even, the canonical class of $X$ is given by $K_X = \rho^*(K_{\ff_{2k}} + L)$ where $K_{\ff_{2k}} = -2 B + (2k - 2) F$ and $2L=[C \cup E]$.
Thus we recover the formula used in the proof of Lemma~\ref{lem:ornor} from $K_X = \rho^*((k- 2)F)$. We therefore conclude that $w_1(\xr) = 0$ hence
$\xr$ is orientable.
\end{proof}

%%%%%%%%%%%%%%%%%%%%%%
\section{Constructions}\label{sec:collapsing}
%%%%%%%%%%%%%%%%%%%%%%

\subsection{Collapsing ovals of trigonal curves via dessins d'enfants}

Let $C$ be a real nonsingular trigonal curve on $\ff_{2k}$ with affine reduced equation

\begin{equation}\label{E:affineequationforC}
x^3+p(u)x+q(u)=0
\end{equation}
where $\mbox{deg} \, p= 4k$ and $\mbox{deg} \, q=6k$.
We will assume that $C$ is generic in the sense that the discriminant $\Delta=4p^3+27q^2$
has only simple roots and $p,q,\Delta$ have distinct roots. The roots of $\Delta$ being simple roots, the curve $C$ is a nonsingular curve
and any root of $\Delta$ corresponds to a vertical line tangent to $C$.

Following Orevkov \cite{Or}, consider the real rational function $\pp^1(\cc) \rightarrow \pp^1(\cc)$
defined by
$$f=\frac{\Delta}{q^2}=\frac{4p^3}{q^2}+27 \;. $$
This a rational function of degree $12k$. The roots of $\Delta$, $p$ and $q$ are mapped to $0,27$ and $\infty$, respectively. These roots as well as the values
$0,27$ and $\infty$ will be called {\it special}.

\begin{rem}

Recall that the functional invariant $J: \pp^1(\cc) \rightarrow \pp^1(\cc)$ associated with the elliptic surface
given by the affine equation $y^2=x^3+p(u)x+q(u)$ is

$$
J=\frac{4p^3}{\Delta} \; .
$$

The real rational functions $f$ and $J$ are related via
$J=\gamma \circ f$, where $\gamma$ is the real automorphism of $\pp^1(\cc)$
sending $0$, $27$, and $\infty$ to $\infty$, $0$, and $1$, respectively. 
\end{rem}
\medskip

\begin{figure}[htbp]
\bigskip

\begin{center}
\includegraphics{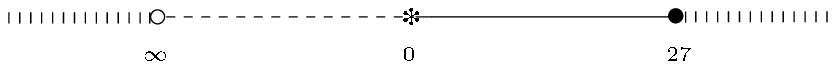}
\caption{Coloring of $\pp^1(\rr)$}
\label{F:coloring}
\end{center}
\bigskip

\end{figure}

Color the real line on the target space in three distinct colors as shown in Figure~\ref{F:coloring} (one color for each interval
delimited by two special values, and special symbols for special values) and define a colored graph on the source space $\pp^1(\cc)$ as the graph $\Gamma=f^{-1}(\pp^1(\rr))$
colored via the pull-back of the previous coloring of $\pp^1(\rr)$. Such a graph is sometimes called a {\it dessin d'enfant}.
The graph $\Gamma$ is invariant under the complex conjugation and contains $\pp^1(\rr)$. Each vertex of $\Gamma$  has even valency and corresponds
to a critical point of $f$ with multiplicity half the valency (if the valency is $2$ then the point is in fact not a critical point of $f$).
Due to the special form of the rational function $f$, the vertices given by the roots of $\Delta$, $p$ and $q$ should have valencies
which are multiple of $2$, $6$ and $4$, respectively. In fact, since we have assumed that $\Delta$ has only simple roots, the roots of $\Delta$
are not critical points of $f$ and have valency $2$.

Hence, to a real trigonal curve $C$ corresponds a colored graph on $\pp^1(\cc)$ verifying the above properties. Conversely, drawing dessins d'enfants allows to
construct real trigonal curves with a given real scheme (see, for example, \cite{Or}, \cite{BeBr}).
Here we will only use a local description of the graph $\Gamma$ that will allow us to collapse independently the ovals.
\smallskip

\begin{prop}\label{P:collapsing}
The ovals of $C$ can be collapsed independently. In other words, if $C$ has real scheme $\langle a \mid b \rangle$,
then for any pair of integers $(a',b')$ such that $0 \leq a' \leq a$ and $0 \leq b' \leq b$ there exists
a real nonsingular trigonal curve $C'$ on $\ff_{2k}$ with real scheme $\langle a' \mid b' \rangle$.
\end{prop}

\begin{figure}[htbp]
\begin{center}
\includegraphics{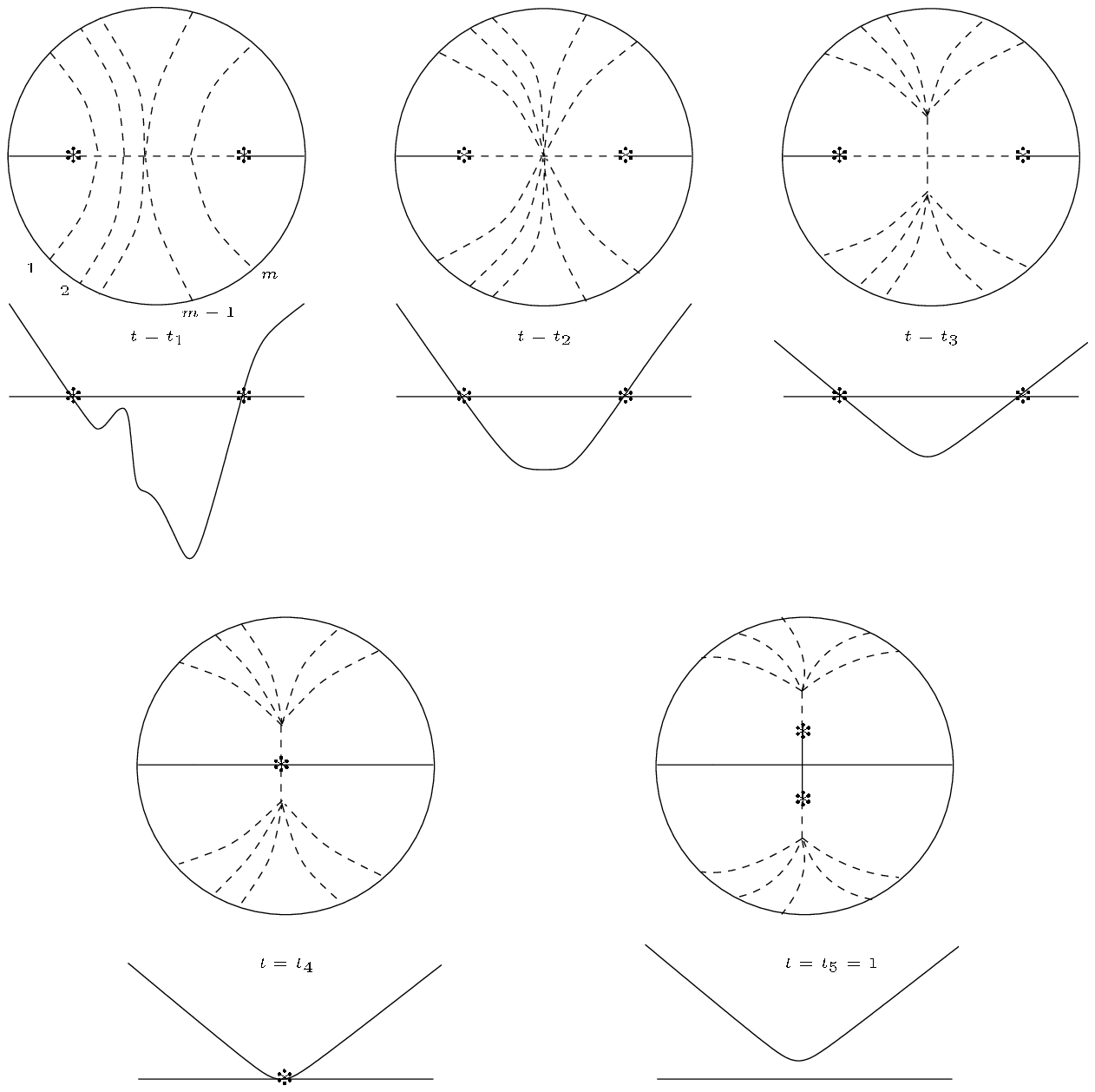}
\caption{}
\label{F:deformation}
\end{center}
\end{figure}

\begin{proof}
Consider one oval of the curve (if the curve has no oval, there is nothing to prove).
Its image under the ruling $\tilde \pi \colon (u,x) \mapsto u$ is a closed interval whose endpoints are roots of $\Delta$.
Denote by $I$ the interior of this interval. For each $u_0 \in I$, the vertical line $u=u_0$ intersects the oval in three
distinct real points. This means that $f$ is negative on $I$. In other words, the interval $I$ is uniquely colored as a part of $\Gamma$
like $(\infty, 0)$ on the target space $\pp^1(\rr)$ (see Figure~\ref{F:coloring}).
Conversely, any such interval corresponds to an oval of the curve.
The critical points of $f$ which belong to $I$ are either roots of $q$ (critical points with value $\infty$),
or non special critical points.
Let $\mathcal U$ be a small neighbourhood of $I \subset \pp^1(\cc)$ which is invariant under the complex conjugation and which
contains only real critical points of $f$.  Let $H$ be one of the two complex conjugate parts of $\mathcal U \setminus \pp^1(\rr)$.
Recall that the valency of a vertex of $\Gamma$ is twice its multiplicity as a critical point of $f$.
The multiplicity of each root of $q$ contained in $I$ should be even for otherwise $f$ would be positive inside $I$.
This means that the number of branches of $\Gamma$ contained in $H$ and which start from a real root of $q$ is odd.
Moreover, since the derivative $f'$ take opposite (non zero) values at the endpoints of $I$, it follows that
the sum of the multiplicities as roots of $f'$ of the non special critical points contained in $I$ should be an odd integer $m$.
We note that this $m$ coincides with the total number of branches of $\Gamma$ contained in $H$ and which start from these non special critical points.

\begin{figure}[htbp]

\begin{center}
\includegraphics{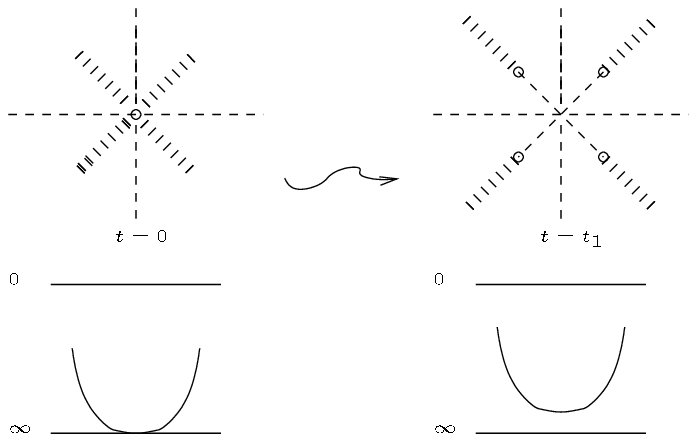}
\caption{Deformation around a real root of $q$ (with multiplicity $4$)}
\label{F:deformationbis}
\end{center}

\end{figure}

We perform a continuous deformation of colored graphs ${(\Gamma_t)}_{t \in [0,1]}$, $\Gamma=\Gamma_0$, which is concentrated in
the neighbourhood $\mathcal U$, as it is shown in Figure~\ref{F:deformationbis} and~\ref{F:deformation} (the corresponding deformation
for the graph of $f_{|\pp^1(\rr)}$ is depicted at the bottom).
In these pictures, horizontal segments are parts of $\pp^1(\rr)$ while the complex conjugation is given by the reflection about $\pp^1(\rr)$.

First we deform $\Gamma$ around each real root of $q$ as depicted in Figure~\ref{F:deformationbis}.
The real root is transformed into pairs of complex conjugate simple roots, while a non special real critical point
(of the same mulitplicity) is created.
Shrinking $\mathcal U$ if necessary, we can then assume that it contains only (real) non special critical points of $f$.
This means that, with the exception of the complement of $I$ in $\pp^1(\rr) \cap \mathcal{U}$, the part of $\Gamma$ contained in $\mathcal U$ is now
entirely colored like $(\infty, 0)$ and the total number of branches contained in $H$ is the odd number $m$. This gives a new graph at time $t=t_1$.
We then pursue the deformation as shown in Figure~\ref{F:deformation} (the disk represents $\mathcal{U}$).

It remains to show that along the deformation any graph $\Gamma_t$ corresponds to a real trigonal curve $C_t$ on $\ff_{2k}$
and that the real scheme of $C_1$ is obtained from that of $C_0$ by collapsing the oval we started with.

From the topological point of view, the map $f: \pp^1(\cc) \rightarrow \pp^1(\cc)$ is a symmetric (by which we mean
that complex conjugate points are mapped to complex conjugate points) branched covering
of degree $12k$ which respects the colorings.
It is not difficult to see that the deformation ${(\Gamma_t)}_{t \in [0,1]}$ comes from a continous
deformation $(f_t)_{t \in [0,1]}$, with $f_0=f$, of maps $f_t$ sharing the same properties and such that $\Gamma_t= f_t^{-1}(\pp^1(\rr))$.
Indeed, there exists a continuous map
$G: \pp^1(\cc) \times [0,1] \rightarrow \pp^1(\cc)$ such that $\Gamma_t=G(\Gamma_0,t)$ for each $t$
and the restriction of $G$ to $\Gamma_0 \times \{0\}$ is given by the identity map.
Moreover, there is a continous family ${(\varphi_t)}_{t \in [0,1]}$ of maps $\varphi_t: \Gamma_t \rightarrow \pp^1(\rr)$
starting at the restriction $\varphi_0$ of $f$ to $\Gamma=\Gamma_0$ and which respect the colorings.
Consider now one connected component $D_0$ of $\pp^1(\cc) \setminus \Gamma_0$ and let $D_t=G(D_0,t)$ be the corresponding
connected component of $\pp^1(\cc) \setminus \Gamma_t$. Denote by $\bar{D}_t$ and $\bar{D}_0$ their closures.
Looking at Figure~\ref{F:deformation} shows that $\bar{D}_t$ and $\bar{D}_0$ are homeomorphic
for each $t$. Moreover, one can choose homeomorphisms $\psi_{D,t}: \bar{D}_t \rightarrow\bar{D}_0 $ coming into a continuous family such that
${\varphi_t}$ and $\varphi_0 \circ {\psi_{D,t}}$ coincide on the boundary of $D_t$. This allows to extend the restriction
of $\varphi_t$ to the boundary of $D_t$ into a continuous map $f_{D,t}: \bar{D}_t \rightarrow f_0({\bar{D}_0})$
by setting $f_{D,t}:= f_0 \circ \psi_{D,t}$. Doing that for each $D_0$ in one given connected component of
$\pp^1(\cc) \setminus \pp^1(\rr)$, and extending the resulting map into a symmetric map, one obtains a continuous family of
symmetric continuous maps $f_t: \pp^1(\cc) \rightarrow \pp^1(\cc)$. It is then easy to see that each $f_t$ has the desired properties.

By Riemann's uniformization Theorem, each map $f_t$ becomes a real rational map of degree $12k$ for the standard complex structure
on the target space and its pull-back by $f_t$ on the source space.
For any $t \in [0,1]$ the (special) points of $\Gamma_t$ where $f_t$ take values $0$,
$\infty$ and $27$ have all a valency which is a multiple of $2$, $4$ and $6$, respectively.
This gives the existence of real polynomials $p_t$ and $q_t$ of degree $4k$ and $6k$, respectively
\footnote{A priori we can only conclude that $\mbox{deg} \, p_t \leq 4k$ and $\mbox{deg} \, q_t \leq 6k$
with at least one equality occuring, but Figure~\ref{F:deformation} shows that
$p_t$ and $p_0=p$ (resp., $q_t$ and $q_0$) have the same number of roots, counted with multiplicities,
hence the same degree.}, such that $f_t=\frac{4p_t^3}{q_t^2}+27=\frac{\Delta_t}{q_t^2}$ where $\Delta_t:=4p_t^3+27q_t^2$.
Let $C_t$ denotes the associated trigonal curve with affine equation
$x^3+p_t(u)x+q_t(u)=0$. For any $t \neq t_4$, the curve $C_t$ is a nonsingular curve since all roots of $\Delta_t$
have valency $2$ in $\Gamma_t$. The curve $C_{t_4}$ has an unique singular point which corresponds
to the real root of $\Delta_{t_4}$ contained in $\mathcal U$ whose valency is $4$.
This singular point is an ordinary double point which is isolated in $\ff_{2k}(\rr)$ since $f_{t_4}$, hence $\Delta_{t_4}$, is positive around it.
To finish, it suffices to note that the number of intervals uniquely colored like $(0,\infty)$ on the target space
decreases by $1$ along the deformation.
\end{proof}

%%%%%%%%%%%%%%%%%%%%%%%%%%%%%%%%%%%%%%%%
\subsection{Combinatorial patchworking of trigonal curves}\label{S:patchworking}

We first recall how the combinatorial patchworking works for the construction of real trigonal curves on $\ff_{2k}$
Let $k$ be any positive integer and denote by $T$ the triangle with vertices
$(0,0)$, $(6k,0)$ and $(0,3)$. The input of the combinatorial patchworking are
\begin{itemize}
\item A triangulation of $T$ whose vertices have integer coordinates,
\item a distribution of signs $\pm 1$ at the vertices of this triangulation.
\end{itemize}

Assume that there exists a convex piecewise-linear function
$\nu : T \to {\rr}$ linear on each triangle of the triangulation,
but not linear on the union of any two triangles. Such a triangulation is called {\it convex} (or sometimes {\it coherent}).

The so-called T-polynomial associated with these data is the polynomial $r_t \in \rr[u,x]$ defined by

$$r_t(u,x)= \sum s(i_1,i_2) t^{\nu(i_1,i_2)} u^{i_1} x^{i_2},$$
where the sum is taken over all the vertices $(i_1,i_2)$ of the triangulation and $s(i_1,i_2)$
is the sign of $(i_1,i_2)$.

Consider now the quadrangle $Q$ with vertices $(\pm 6k,0)$ and $(0, \pm 3)$.
Extend the triangulation of $T$ to a triangulation of $Q$ which is symmetric with respect to the coordinate axes.
Extend also the initial sign distribution to a sign distribution at the vertices of the triangulation of $Q$
following the rule
$$s(\epsilon_1 i_1,\epsilon_2 i_2):=\epsilon_1^{i_1} \cdot \epsilon_2^{i_2} \cdot s(i_1,i_2)$$
for any $(\epsilon_1,\epsilon_2) \in \{\pm 1\}^2$.

We are now able to perform the combinatorial patchworking.
For any triangle of the triangulation of $Q$ whose vertices do not all have the same sign,
select the segment joigning the middle points of the two edges whose endpoints have different signs.
Denote by $\mathcal C$ the piecewise linear curve in $Q$ obtained by taking the union of
all these selected segments. We will use the following version of the combinatorial patchworking theorem, see ~\cite{Vi83, Vi84, Vi}.

\begin{prop}\label{P:patchwork} There exists $t_0>0$ such that for any positive $t <t_0$ the polynomial $r_t$
defines a real nonsingular trigonal curve $C$. Moreover, there exists an homeomorphism
from ${\rr}^2$ to the interior of $Q$ sending the affine curve defined by $r_t$ to $\mathcal C$.
\end{prop}

\begin{figure}
\bigskip

\begin{center}
\includegraphics{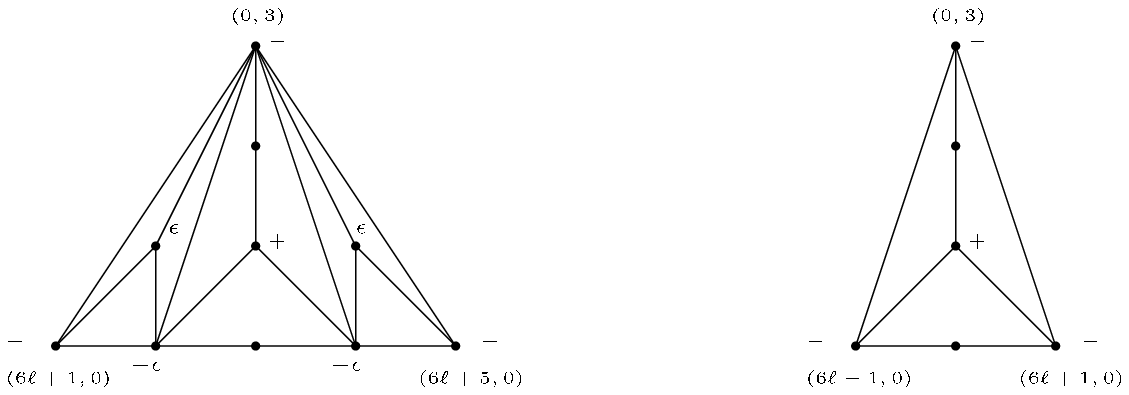}
\caption{Triangulation and signs For $T_{\ell}$ and $t_{\ell}$.}
\label{F:triang}
\end{center}
\begin{center}
{\scriptsize Points with integer coordinates are represented by black points.}
\end{center}
\bigskip

\end{figure}

We recall the Harnack's theorem for real curves, which is a particular case of Proposition~\ref{prop:smith}. If $C$ is
a real nonsingular algebraic curve, then the number of connected components of $C(\rr)$ is no more than
$g(C)+1$, where $g(C)$ is the genus of $C$. If $C(\rr)$ has $g(C)+1$ connected components, then $C$ is called an $M$-curve,
while it is called an $(M-d)$-curve if $C(\rr)$ has $g(C)+1-d$ connected components.

\begin{prop}\label{P:Mcurve}
\begin{enumerate}
\item For any integer $\lambda$ such that $0 \leq \lambda \leq k$, there exists a real nonsingular
trigonal $M$-curve $C$ with real scheme $\langle k-1+4\lambda \mid 5k-1-4\lambda \rangle$.
\item For any integer $\lambda$ such that $0 \leq \lambda \leq k-1$, there exists a real nonsingular
trigonal $(M-2)$-curve $C$ on $\ff_{2k}$ with real scheme $\langle k+4\lambda \mid 5k-4-4\lambda \rangle$.
\end{enumerate}
\end{prop}

\begin{figure}
\bigskip

\begin{center}
\includegraphics{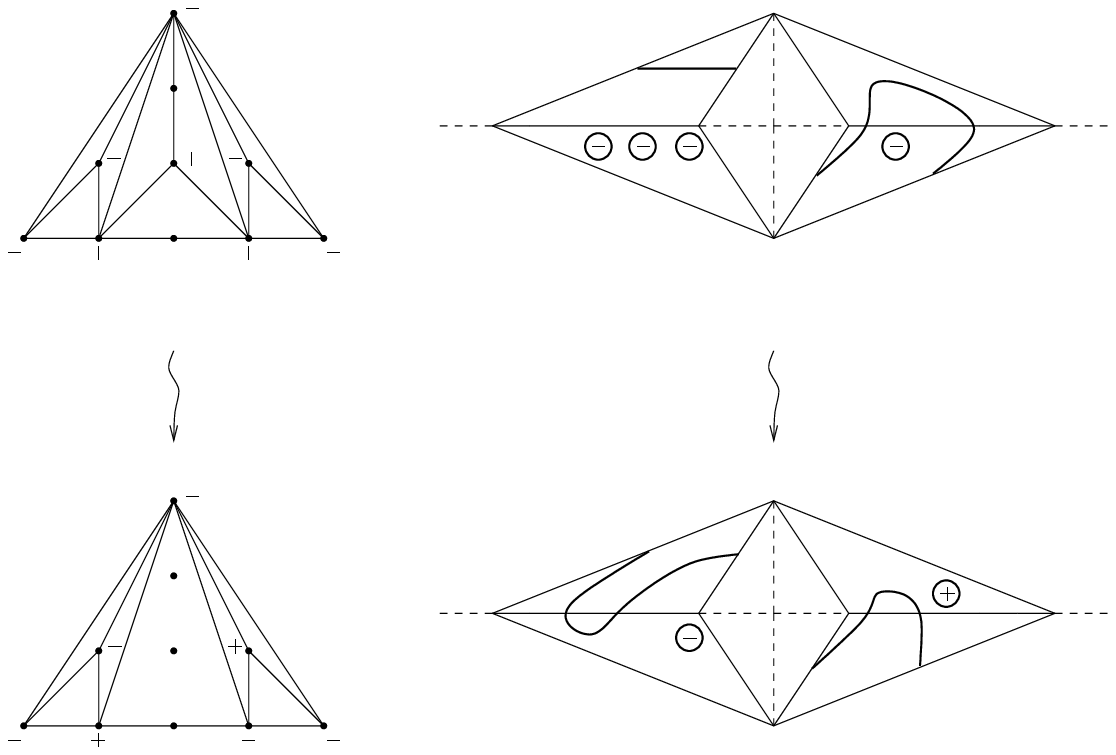}\caption{Modification of the triangulation and sign distribution}
\label{F:modiftriang}
\end{center}
\bigskip

\end{figure}

\begin{proof}
Use the segments $[(0,3),(6 \ell \pm 1)]$ in order to triangulate $T$ into the two triangles
$[(0,0), (1,0), (0,3)]$, $[(6k-1,0), (6k,0), (0,3)]$ and the triangles

$$
\begin{array}{llll}
T_{\ell} & := & [(6\ell+1, 0), (6\ell+5,0), (0,3)],  & \quad \ell=0,\dotsc, k-1 \\
 & & & \\
t_{\ell} & := & [(6\ell-1,0), (6\ell+1,0), (0,3)], & \quad \ell=1,\dotsc, k-1.
\end{array}
$$
Then triangulate each $T_{\ell}$ and $t_{\ell}$ as depicted in Figure~\ref{F:triang}.
It is easy to see that the resulting triangulation of $T$ is convex, see Remark~\ref{R:convex}.

Let us now choose a sign distribution at the vertices of this triangulation.
Put the sign $-1$ to the points $(0,0)$, $(6k,0)$ and take the sign distribution
shown in  Figure~\ref{F:triang} for the vertices of the triangulation which are contained in triangles $t_{\ell}$.
For the vertices of the triangulation which are contained in triangles $T_{\ell}$, take the sign distribution shown in Figure~\ref{F:triang} with $\epsilon=+$
for $\lambda$ triangles among $T_0, \dotsc, T_{k-1}$ and with $\epsilon=-$ for the remaining $k-\lambda$ triangles.

Let $r_t \in \rr[u,x]$ be a T-polynomial associated with these data.
Let $C$ be the trigonal curve on $\ff_{2k}$ defined by $r_t$ for $t>0$ sufficiently small and let ${\mathcal C} \subset Q$ be the piecewise linear curve
constructed via the combinatorial patchworking from the previous triangulation and sign distribution.
Proposition~\ref{P:patchwork} implies that $C$ is a nonsingular curve whose desired properties can be read off the picture of ${\mathcal C} \subset Q$.

One can easily check that the union of the four symmetric copies of the interior of each triangle $t_{\ell}$ contains exactly two ovals of $\mathcal C$,
one of them surrounding a vertex with the sign $+$ while the other surrounds a vertex with the sign $-$.
It is also easy to see that the union of the four symmetric copies of the interior of each triangle $T_{\ell}$ contains four ovals of $\mathcal C$,
all of them surrounding a vertex with the sign $\epsilon$ (see Figure~\ref{F:modiftriang} for $\epsilon=-1$).
This gives $2(k-1)+4k = 6k-2$ ovals of $C(\rr)$.
Note that $\mathcal C$ has one connected component more which intersects the segment with vertices $(0,\pm 3)$. This gives an account of $6k-1$ connected components of
$C(\rr)$. But $6k-2$ is the genus of the complex curve $C$. Hence $C(\rr)$ has exactly $6k-1$ connected components by Harnack's theorem, that is $C$ is an $M$-curve,
and the last component we have obtained is the pseudo line.
Finally,  it remains to compute the real scheme $\langle a \mid b \rangle$ of $C$. Clearly, up to permutation of $a$ and $b$,
the number $a$ is the number of ovals of $\mathcal C$ which surround a vertex with the sign $+$, while $b$ is the number of ovals of $\mathcal C$
which surround a vertex with the sign $-$. Since we have choosen the sign $\epsilon=+$ for $\lambda$ triangles among $T_0,\dotsc, T_{k-1}$,
we obtain $a=k-1+4\lambda$ and $b=k-1+4(k-\lambda)=5k-1-4\lambda$.

In order to prove the statement about $(M-2)$-curves,
we modify slightly the previous triangulation and sign distribution.
Consider one triangle among those $k-\lambda >0$ triangles $T_{\ell}$ for which the sign $\epsilon=-1$ has been choosen
and modify its triangulation and corresponding sign distribution as depicted in Figure~\ref{F:modiftriang}.
The corresponding modification of real scheme is shown in Figure~\ref{F:modiftriang}.
It follows that this modification increases by $1$ the number $a$ and decreases by $3$ the number $b$.
\end{proof}

\begin{rem}\label{R:convex}
A classical argument explained below shows that all triangulations we used are convex.
Note that this is in fact unnecessary due to~\cite{BeBr} which proves that the convexity hypothesis can be dropped
in the combinatorial patchworking theorem when starting from a triangulation of $[(0,0), (3n,0),(0,3)]$.
Take a triangulation of $T$ by segments starting from $(0,3)$ and ending to points of $[(0,0), (6k,0)]$. Such a triangulation is evidently convex.
Now subdivide one its triangles into three triangles using an interior point. Let $\nu$ be a function which certifies the convexity of
the starting triangulation. We can construct a function $\nu'$ certifying the convexity of the refined triangulation in the following way:
let $\nu'$ coincide with $\nu$ at the vertices of the starting triangulation, take the same value than $\nu$ but diminished by $\iota$ with $0< \iota <<1$ at the new vertex,
and be linear on each triangle of the new triangulation.
The triangulations used in this paper are all obtained by successive refinements of the previous type, and are thus convex.
\end{rem}

%%%%%%%%%%%%%%%%%%%%%%%%%%%%%%%%%%%%%
\subsection{Proof of Theorem~\ref{theo:classif}}\label{subsection:finalproof}

We organize all the results obtained so far to give a proof of Theorem~\ref{theo:classif}.
The fact that the topological types of real jacobian elliptic surfaces of given holomorphic characteristic $k$
belong to the list in Theorem~\ref{theo:classif} has been proven in Section~\ref{sec:restrict}.

By Lemma~\ref{L:curvetosurf}, the extremal topological type $S_1 \sqcup S_1$ if $k$ is even, or $V_2 \sqcup V_2$ if $k$ is odd,
is obtained via a double covering starting with a trigonal curve constructed in Example~\ref{E:specialextremal}.

Applying Lemma~\ref{L:curvetosurf} to the curves in parts 1) and 2) of Proposition~\ref{P:Mcurve}
produces via a double covering real jacobian elliptic surfaces
whose topological types form the list $(1)$ and $(2)$ of Theorem~\ref{theo:classif}, respectively.
Now, Proposition~\ref{P:collapsing} applies to each of these curves, which in turn implies
via Lemma~\ref{L:curvetosurf} that all desired topological types can be realized by real jacobian elliptic surfaces
with holomorphic Euler characteristic $k$.

%%%%%%%%%%%%%%%%%%%%%%%%%%%%%%%%%%%%%%%%

\end{document}